\documentclass[12pt]{article}

\usepackage{amssymb}
\usepackage{amsmath}
\usepackage{cite}
\usepackage[arrow, matrix, curve]{xy}

\begin{document}

\renewcommand{\citeleft}{{\rm [}}
\renewcommand{\citeright}{{\rm ]}}
\renewcommand{\citepunct}{{\rm,\ }}
\renewcommand{\citemid}{{\rm,\ }}

\newcounter{abschnitt}
\newtheorem{satz}{Theorem}
\newtheorem{theorem}{Theorem}[abschnitt]
\newtheorem{koro}[theorem]{Corollary}
\newtheorem{prop}[theorem]{Proposition}
\newtheorem{lem}[theorem]{Lemma}
\newtheorem{conj}[theorem]{Conjecture}
\newtheorem{probl}[theorem]{Problem}

\newcounter{saveeqn}
\newcommand{\alpheqn}{\setcounter{saveeqn}{\value{abschnitt}}
\renewcommand{\theequation}{\mbox{\arabic{saveeqn}.\arabic{equation}}}}
\newcommand{\reseteqn}{\setcounter{equation}{0}
\renewcommand{\theequation}{\arabic{equation}}}

\hyphenation{convex} \hyphenation{bodies}

\sloppy

\phantom{a}

\vspace{-2.2cm}

\begin{center}
 \begin{LARGE} { The Steiner Formula for Minkowski Valuations} \\[0.7cm] \end{LARGE}

\begin{large} Lukas Parapatits and Franz E.
Schuster \end{large}
\end{center}

\vspace{-1.1cm}

\begin{quote}
\footnotesize{ \vskip 1truecm\noindent {\bf Abstract.} A Steiner
type formula for continuous translation invariant Minkowski
valuations is established. In combination with a recent result on
the symmetry of rigid motion invariant homogeneous bivaluations,
this new Steiner type formula is used to obtain a family of
Brunn--Minkowski type inequalities for rigid motion intertwining
Minkowski valuations. }
\end{quote}

\vspace{0.6cm}

\centerline{\large{\bf{ \setcounter{abschnitt}{1}
\arabic{abschnitt}. Introduction}}} \alpheqn

\vspace{0.6cm}

The famous Steiner formula, dating back to the 19th century,
expresses the volume of the parallel set of a convex body $K$ at
distance $r \geq 0$ as a polynomial in $r$. Up to constants
(depending on the dimension of the ambient space), the
coefficients of this polynomial are the intrinsic volumes of $K$.
Steiner's formula is among the most influential results of the
early days of convex geometry. Its ramifications and many
applications can be found, even today, in several mathematical
areas such as differential geometry (starting from Weyl's tube
formula \textbf{\cite{weyl39}}; see e.g.\ \textbf{\cite{fu85,
gray90}} for more recent results), geometric measure theory
(going back to Federer's seminal work on curvature measures
\textbf{\cite{federer59}}; see also \textbf{\cite{fu94a, fu94b,
schneider79, ratzaehl03}}), convex and stochastic geometry (see
e.g.\ \textbf{\cite{hulawe04, schneider93, schnweil}}), geometric
functional analysis (see \textbf{\cite{colesanti, colhug00}}),
and recently also in algebraic geometry (see
\textbf{\cite{henher08, teissier82}}).

In Euclidean space $\mathbb{R}^n$, the parallel set of $K$ at
distance $r$ is equal to the sum of $K$ and a Euclidean ball of
radius $r$. A fundamental extension of the classical Steiner
formula is Minkowski's theorem on the polynomial expansion
\linebreak of the volume of a Minkowski sum of several convex
bodies, leading to the theory of mixed volumes (see e.g.\
\textbf{\cite{schneider93}}). More recently, McMullen
\textbf{\cite{McMullen77}} (and later, independently, Meier
\textbf{\cite{meier77}} and Spiegel \textbf{\cite{spiegel78}})
established the existence \linebreak of a similar polynomial
expansion for functions on convex bodies which are considerably
more general than the ordinary volume, namely continuous
translation invariant (real valued) {\it valuations}.

The origins of the notion of valuation (see Section 2 for precise
definition) can be traced back to Dehn's solution of Hilbert's
Third Problem. However, \linebreak the starting point for a
systematic investigation of general valuations was Hadwiger's
\textbf{\cite{hadwiger51}} fundamental characterization of the
linear combinations of intrinsic volumes as the continuous
valuations that are rigid motion invariant (see
\textbf{\cite{Alesker99, Alesker01, Bernig09, centro}} for recent
important variants). McMullen's \textbf{\cite{McMullen77}} deep
result on the polynomial expansion of translation invariant
valuations is among the seminal contributions to the structure
theory of the space of translation invariant valuations which has
been rapidly evolving over the last decade (see
\textbf{\cite{Alesker01, Alesker03, AlBern, Bernig07b, fu06}}).
These recent structural insights in turn provided the means for a
fuller understanding of the integral geometry of groups acting
transitively on the sphere (see e.g.\ \textbf{\cite{Alesker03,
Bernig09, bernig10, bernigfu10}} and the survey
\textbf{\cite{bernig10a}}).

While classical results on valuations were mainly concerned with
real and tensor valued valuations, a very recent development
explores the strong connections between {\it convex body valued}
valuations and isoperimetric and related inequalities (see
\textbf{\cite{ABS2011, habschu09, Ludwig10, Schu09}}). This new
line of research has its roots in the work of Ludwig
\textbf{\cite{ludwig02, Ludwig:matrix, Ludwig:Minkowski,
Ludwig06}} who first obtained classifications \linebreak of
convex and star body valued valuations which are compatible with
linear transformations (see also \textbf{\cite{haberl08,
haberl09, haberl10, hablud06, Ludwig10a, wannerer10, SchuWan11}}).
In this area, it is a major \linebreak open problem whether a
polynomial expansion of translation invariant convex \linebreak
body valued valuations is also possible (see Section 2 for
details).

In this article we establish a Steiner type formula for
continuous translation invariant {\it Minkowski valuations}
(i.e.\ valuations taking values in the topological semigroup of
convex bodies endowed with Minkowski addition). In fact, we
obtain a more general polynomial expansion formula for translation
invariant Minkowski valuations when the arguments are Minkowski
sums of zonoids. This follows in part from a connection between
Minkowski \linebreak valuations and positive scalar valuations.
Our new Steiner type formula gives rise to a Lefschetz operator on
Minkowski valuations which we use \linebreak together with a
recent result on the symmetry of rigid motion intertwining
\linebreak homogeneous bivaluations \textbf{\cite{ABS2011}} to
obtain a family of Brunn--Minkowski type inequalities for
intrinsic volumes of rigid motion intertwining Minkowski
\linebreak valuations. These new inequalities generalize a number
of previous partial results \textbf{\cite{ABS2011, lutwak93,
Schu06b, Schu09}}.

\vspace{1cm}

\centerline{\large{\bf{ \setcounter{abschnitt}{2}
\arabic{abschnitt}. Statement of principal results}}} \alpheqn

\vspace{0.6cm}

The setting for this article is $n$-dimensional Euclidean space
$\mathbb{R}^n$ with $n \geq 3$. We denote by $\mathcal{K}^n$ the
space of convex bodies in $\mathbb{R}^n$ endowed with the
Hausdorff metric. A function $\varphi$ defined on $\mathcal{K}^n$
and taking values in an abelian semigroup is called a {\it
valuation} if
\[\varphi(K) + \varphi(L) = \varphi(K \cup L) + \varphi(K \cap L)\]
whenever $K \cup L \in \mathcal{K}^n$. A valuation $\varphi$ is
said to be {\it translation invariant} if $\varphi(K + x) =
\varphi(K)$ for all $x \in \mathbb{R}^n$ and $K \in
\mathcal{K}^n$.

The most familiar real valued valuation is, of course, the
ordinary volume $V_n$. In fact, the valuation property of volume
carries over to a series of basic functions which are derived
from it: By a classical result of Minkowski, the volume of a
Minkowski (or vector) linear combination $\lambda_1K_1 + \cdots +
\lambda_mK_m$ of convex bodies $K_1, \ldots, K_m \in
\mathcal{K}^n$ with real coefficients $\lambda_1, \ldots,
\lambda_m \geq 0$ can be expressed as a homogeneous polynomial of
degree $n$,
\begin{equation} \label{mixed}
V_n(\lambda_1K_1 + \cdots +\lambda_m K_m)=\sum
\limits_{j_1,\ldots, j_n=1}^m
V(K_{j_1},\ldots,K_{j_n})\lambda_{j_1}\cdots\lambda_{j_n},
\end{equation}
where the coefficients $V(K_{j_1},\ldots,K_{j_n})$, called {\it
mixed volumes} of $K_{j_1}, \ldots, K_{j_n}$, are symmetric in the
indices and depend only on $K_{j_1}, \ldots, K_{j_n}$. Now, if
\linebreak $i \in \{1, \ldots, n\}$ and an $(n-i)$-tuple $L_1,
\ldots, L_{n-i}$ of convex bodies is fixed, then the function
$\phi: \mathcal{K}^n \rightarrow \mathbb{R}$, defined by
$\phi(K)=V(K,\ldots,K,L_1,\ldots,L_{n-i})$, is a continuous
translation invariant valuation (see e.g.\
\textbf{\cite{schneider93}}).

\vspace{0.1cm}

In a highly influential article, Alesker \textbf{\cite{Alesker01}}
showed (thereby confirming a conjecture by McMullen) that in fact
every continuous translation invariant real valued valuation is a
limit of linear combinations of mixed volumes. One of the crucial
ingredients in the proof of Alesker's landmark result is the
following significant generalization of the polynomial expansion
(\ref{mixed}):

\vspace{0.1cm}

\begin{satz} \label{polyexp} {\bf (McMullen \cite{McMullen77})}
Let $X$ be a topological vector space. Suppose that $\varphi:
\mathcal{K}^n \rightarrow X$ is a continuous translation
invariant valuation and let $K_1, \ldots, K_m \in \mathcal{K}^n$.
Then
\[\varphi(\lambda_1K_1 + \cdots + \lambda_mK_m), \qquad \lambda_1, \ldots, \lambda_m \geq 0,\]
can be expressed as a polynomial in $\lambda_1, \ldots, \lambda_m$
of total degree at most $n$. Moreover, for each
$(i_1,\ldots,i_m)$, the coefficient of
$\lambda_1^{i_1}\cdots\lambda_m^{i_m}$ is a continuous
translation invariant and homogeneous valuation of degree $i_j$
in $K_j$.
\end{satz}

\vspace{0.1cm}

As a special case of Theorem \ref{polyexp}, we note the following
extension of the classical Steiner formula for volume (see
Section 5): If $K \in \mathcal{K}^n$, then for every $r \geq 0$,
\begin{equation} \label{mcmullenstein}
\varphi(K + rB^n) = \sum_{j=0}^n r^{n-j}\varphi^{(j)}(K),
\end{equation}
where the coefficient functions $\varphi^{(j)}: \mathcal{K}^n
\rightarrow X$, $0 \leq j \leq n$, defined by
(\ref{mcmullenstein}), are continuous translation invariant
valuations. Clearly, $\varphi^{(n)}=\varphi$.

\vspace{0.3cm}

\noindent {\bf Definition} \emph{A map $\Phi: \mathcal{K}^n
\rightarrow \mathcal{K}^n$ is called a Minkowski valuation if
\[\Phi(K) + \Phi(L) = \Phi(K \cup L) + \Phi(K \cap L),  \]
whenever $K, L, K \cup L \in \mathcal{K}^n$ and addition on
$\mathcal{K}^n$ is Minkowski addition.}

\vspace{0.3cm}

While first results on Minkowski valuations were obtained in the
1970s by Schneider \textbf{\cite{schneider74}}, they have become
the focus of increased interest (and acquired their name) more
recently through the work of Ludwig \textbf{\cite{ludwig02,
Ludwig:Minkowski}}. It was shown there that such central notions
like projection, centroid and difference body operators can be
characterized as unique Minkowski valuations compatible with
affine transformations of $\mathbb{R}^n$ (see
\textbf{\cite{haberl09, hablud06, Ludwig10a, SchuWan11,
wannerer10}} for related results).

Since the space of convex bodies $\mathcal{K}^n$ does not carry a
linear structure, it is an important open problem (cf.
\textbf{\cite{Schu09}}) whether the Steiner type formula
(\ref{mcmullenstein}), or even Theorem \ref{polyexp}, also hold
for continuous translation invariant Minkowski valuations. As our
main result we establish an affirmative answer to the first
question:

\begin{satz}\label{main1}  Suppose that $\Phi: \mathcal{K}^n \rightarrow \mathcal{K}^n$ is a continuous translation invariant
Minkowski valuation and let $K \in \mathcal{K}^n$. Then $\Phi(K +
rB^n)$, $r \geq 0$, can be expressed as a polynomial in $r$ of
degree at most $n$ whose coefficients are convex bodies, say
\begin{equation} \label{minkstein}
 \Phi(K + rB^n) = \sum_{j=0}^n r^{n-j} \Phi^{(j)}(K).
\end{equation}
Moreover, the maps $\Phi^{(j)}: \mathcal{K}^n \rightarrow
\mathcal{K}^n$, $0 \leq j \leq n$, defined by (\ref{minkstein}),
are also continuous translation invariant Minkowski valuations.
\end{satz}

The proof of Theorem \ref{main1} makes critical use of an
embedding by Klain \textbf{\cite{klain00}} of translation
invariant continuous {\it even} (real valued) valuations in the
space of continuous functions on the Grassmannian. In fact, our
proof yields a stronger result than Theorem \ref{main1}, see
Corollary \ref{minksteinzon}, where the Euclidean unit ball $B^n$
in (\ref{minkstein}) can be replaced by an arbitrary zonoid
(i.e.\ a Hausdorff limit \pagebreak

\noindent of finite Minkowski sums of line segments). Moreover, in
Theorem \ref{minkpoly} we obtain a polynomial expansion formula
for continuous translation invariant Minkowski valuations when
the summands are zonoids. Very recently, during the review
process of this article, Wannerer and the first-named author
\textbf{\cite{papwann12}} showed that a polynomial expansion
(analogous to Theorem~\ref{polyexp}) of continuous translation
invariant Minkowski valuations is in general {\it not} possible.

A special case of Theorem 2 was previously obtained by the second-named
author \textbf{\cite{Schu06a}}, when the Minkowski valuation
$\Phi$ is in addition \emph{$\mathrm{SO}(n)$ equivariant} and has
\emph{degree $n-1$}, i.e.\ $\Phi(\vartheta K)=\vartheta \Phi(K)$
and $\Phi(\lambda K)=\lambda^{n-1}\Phi(K)$ for every $K \in
\mathcal{K}^n$, $\vartheta \in \mathrm{SO}(n)$ and real $\lambda
> 0$. As an application of this particular case of Theorem 2, an array of geometric inequalities for the
intrinsic volumes $V_i$ of the {\it derived} Minkowski valuations
$\Phi^{(j)}$ (of degree $j - 1$) was obtained in
\textbf{\cite{Schu06b}}. In particular, the following
Brunn--Minkowski type inequality was established: If $K, L \in
\mathcal{K}^n$ and $3 \leq j \leq n$, $1 \leq i \leq n$, then
\begin{equation} \label{genbmbmind}
V_i(\Phi^{(j)}(K+L))^{1/i(j-1)} \geq
V_i(\Phi^{(j)}(K))^{1/i(j-1)}+V_i(\Phi^{(j)}(L))^{1/i(j-1)}.
\end{equation}
It was also shown in \textbf{\cite{Schu06b}} that if $\Phi$ is
{\it non-trivial}, i.e.\ it does not map every convex body to the
origin, equality holds in (\ref{genbmbmind}) for convex bodies
$K$ and $L$ \linebreak with non-empty interior if and only if they
are homothetic.

The family of inequalities (\ref{genbmbmind}) extended at the
same time previously established inequalities for projection
bodies by Lutwak \textbf{\cite{lutwak93}} and the famous classical
Brunn--Minkowski inequalities for the intrinsic volumes (see
e.g.\ \textbf{\cite{schneider93}} and the excellent survey
\textbf{\cite{gardner02}}). We conjecture that inequality
(\ref{genbmbmind}) holds in fact for all continuous translation
invariant and $\mathrm{SO}(n)$ equivariant Minkowski valuations
of a given arbitrary degree $j \in \{2, \ldots, n - 1\}$.

Recently, refining the techniques from the seminal work of Lutwak
\textbf{\cite{lutwak93}}, this conjecture was confirmed in the
case $i = j + 1$, first for even valuations in
\textbf{\cite{Schu09}} and subsequently for general valuations in
\textbf{\cite{ABS2011}}. As an application of Theorem 2, we
extend these results to the case $1 \leq i \leq j + 1$.

\begin{satz} \label{main2} Suppose that $\Phi_j: \mathcal{K}^n \rightarrow
\mathcal{K}^n$ is a non-trivial continuous trans\-lation invariant
and $\mathrm{SO}(n)$ equivariant  Minkowski valuation of a given
degree \linebreak $j \in \{2, \ldots, n - 1\}$. If $K, L \in
\mathcal{K}^n$ and $1 \leq i \leq j + 1$, then
\[V_i(\Phi_j(K + L))^{1/ij} \geq V_i(\Phi_j(K))^{1/ij}+V_i(\Phi_j(L))^{1/ij}.  \]
If $K$ and $L$ are of class $C^2_+$, then equality holds if and
only if $K$ and $L$ are homothetic.
\end{satz}

The proof of Theorem \ref{main2} also uses a recent result on the
symmetry of rigid motion invariant homogeneous bivaluations which
we describe in Section 6. For a discussion of the smoothness
assumption, we refer to Section 7.

\vspace{1cm}

\centerline{\large{\bf{ \setcounter{abschnitt}{3}
\arabic{abschnitt}. Background material for the proof of Theorem 2
}}}

\reseteqn \alpheqn \setcounter{theorem}{0}

\vspace{0.6cm}

In this section we first recall some basic facts about convex
bodies and, in particular, zonoids (see, e.g.\
\textbf{\cite{schneider93}}). Furthermore, we collect results on
translation invariant (mostly real valued) valuations needed in
subsequent sections. In particular, we recall an important
embedding of Klain \textbf{\cite{klain00}} of even translation
invariant continuous valuations in the space of continuous
functions on the Grassmannian.

A convex body $K \in \mathcal{K}^n$ is uniquely determined by the
values of its support function $h(K,x)=\max\{x \cdot y: y \in
K\}$, $x \in \mathbb{R}^n$. Clearly, $h(K,\cdot)$ is positively
homogeneous of degree one and subadditive for every $K \in
\mathcal{K}^n$. Conversely, every function with these properties
is the support function of a convex body.

A Minkowski sum of finitely many line segments is called a
zonotope. A convex body that can be approximated, in the
Hausdorff metric, by a sequence of zonotopes is called a {\it
zonoid}. Over the past four decades it has become apparent that
zonoids arise naturally in several different contexts (see e.g.\
\textbf{\cite[\textnormal{Chapter 3.5}]{schneider93}} and the
references therein). It is not hard to show that a convex body $K
\in \mathcal{K}^n$ is an origin-centered zonoid if and only if
its support function can be represented in the form
\[h(K,x)=\int_{S^{n-1}} |x \cdot u|\,d\mu_K(u), \qquad x \in \mathbb{R}^n,  \]
with some even (non-negative) measure $\mu_K$ on $S^{n-1}$. In
this case, the measure $\mu_K$ is unique and is called the {\it
generating measure} of $K$.

\vspace{0.2cm}

We denote by $\mathbf{Val}$ the vector space of continuous
translation invariant {\it real valued} valuations and we use
$\mathbf{Val}_i$ to denote its subspace of all valuations of
degree $i$. Recall that a map $\varphi$ from $\mathcal{K}^n$ to
$\mathbb{R}$ (or $\mathcal{K}^n$) is said to have degree $i$
\linebreak if $\varphi(\lambda K)=\lambda^i \varphi(K)$ for every
$K \in \mathcal{K}^n$ and $\lambda > 0$.
 A valuation $\varphi \in \mathbf{Val}$ is said to be even
(resp.\ odd) if $\varphi(-K)=(-1)^{\varepsilon}\varphi(K)$ with
$\varepsilon = 0$ (resp.\ $\varepsilon = 1$) for every $K \in
\mathcal{K}^n$. We write $\mathbf{Val}_i^+ \subseteq
\mathbf{Val}_i$ for the subspace of even valuations of degree $i$
and $\mathbf{Val}_i^-$ for the subspace of odd valuations of
degree $i$.

\pagebreak

From the important special case $m = 1$ of Theorem \ref{polyexp},
we deduce that if $\varphi \in \mathbf{Val}$, then there exist
(unique) $\varphi_i \in \mathbf{Val}_i$, $0 \leq i \leq n$, such
that
\begin{equation} \label{homdecomp1}
\varphi(\lambda K) = \varphi_0(K) + \lambda \varphi_1(K) + \cdots + \lambda^n\varphi_n(K)
\end{equation}
for every $K \in \mathcal{K}^n$ and $\lambda > 0$. In fact, a
simple inductive argument, shows that (\ref{homdecomp1}) is
equivalent with Theorem \ref{polyexp}. Since, clearly, every real
valued valuation is the sum of an even and an odd valuation, we
immediately obtain the following corollary, known as {\it
McMullen's decomposition} of the space $\mathbf{Val}$:

\begin{koro} \label{mcmullen}
\[\mathbf{Val} = \bigoplus \limits_{i=0}^n \left(\mathbf{Val}_i^+ \oplus \mathbf{Val}_i^-\right).   \]
\end{koro}

It is easy to show that the space $\mathbf{Val}_0$ is
one-dimensional and is spanned by the Euler characteristic $V_0$.
The analogous non-trivial statement for $\mathbf{Val}_n$ was
proved by Hadwiger \textbf{\cite[\textnormal{p.\
79}]{hadwiger51}}:

\begin{lem} \label{hadvol}  If $\varphi \in \mathbf{Val}_n$, then $\varphi$ is a
multiple of the ordinary volume $V_n$.
\end{lem}

Assume now that $\varphi \in \mathbf{Val}_i$ with $1 \leq i \leq
n - 1$. If $K_1, \ldots, K_m \in \mathcal{K}^n$ and $\lambda_1,
\ldots, \lambda_m > 0$, then, by Theorem \ref{polyexp},
\[\varphi(\lambda_1K_1 + \cdots +\lambda_m K_m)=\sum
\limits_{j_1,\ldots, j_i=1}^m
\varphi(K_{j_1},\ldots,K_{j_i})\lambda_{j_1}\cdots\lambda_{j_i},
\] where the coefficients are symmetric in the indices and depend
only on $K_{j_1}, \ldots, K_{j_i}$. Moreover, the coefficient of
$\lambda_1^{i_1}\cdots \lambda_m^{i_m}$, where $i_1 + \cdots + i_m
= i$, is a continuous translation invariant valuation of degree
$i_j$ in $K_j$, called a \linebreak {\it mixed valuation} derived
from $\varphi$. Clearly, we have $\varphi(K,\ldots,K)=\varphi(K)$.

\vspace{0.1cm}

We now turn to Minkowski valuations. Let $\mathbf{MVal}$ denote
the set of \linebreak continuous translation invariant Minkowski
valuations, and write $\mathbf{MVal}_i^{\pm}$ for its subset of
all even/odd Minkowski valuations of degree $i$.

From Lemma \ref{hadvol} and the special case $m = 1$ of Theorem
\ref{polyexp}, applied to valuations with values in the vector
space $C(S^{n-1})$ of continuous functions on $S^{n-1}$, one can
deduce the following decomposition result (cf.\
\textbf{\cite[\textnormal{p.\ 12}]{schnschu}}):

\begin{lem} \label{minkdecom} If $\Phi \in \mathbf{MVal}$, then for every $K \in \mathcal{K}^n$, there exist convex bodies $L_0, L_n
 \in \mathcal{K}^n$ such that
\begin{equation} \label{suppdecomp}
h(\Phi(K),\cdot) = h(L_0,\cdot) + \sum_{i=1}^{n-1} g_i(K,\cdot) +
V(K)h(L_n,\cdot),
\end{equation}
where, for each $i \in \{1, \ldots, n - 1\}$:
\begin{enumerate}
\item[(i)] The function $g_i(K,\cdot)$ is a difference
of support functions.
\item[(ii)] The map $K \mapsto g_i(K,\cdot)$ is a continuous translation
invariant valuation of degree $i$.
\end{enumerate}
\end{lem}

The natural question whether for every $K \in \mathcal{K}^n$, each
function $g_i(K,\cdot)$ is the support function of a convex body
is equivalent to the following problem.

\begin{probl} \label{mainconj} Let $\Phi \in \mathbf{MVal}$ and $K \in
\mathcal{K}^n$. Are there convex bodies $L_0, \Phi_1(K), \ldots,
\Phi_{n-1}(K), L_n \in \mathcal{K}^n$ such that
\begin{equation} \label{conjminkdecomp}
\Phi(\lambda K) = L_0 + \lambda \Phi_1(K) + \cdots +
\lambda^{n-1}\Phi_{n-1}(K) + \lambda^nV(K)L_n
\end{equation}
for every $\lambda > 0$?
\end{probl}

During the review process of this article, Wannerer and the
first-named author \textbf{\cite{papwann12}} showed that the
answer to Problem \ref{mainconj} is in general negativ. However,
in the next section, we show that (\ref{conjminkdecomp}) holds
for every $\Phi \in \mathbf{MVal}$ and $\lambda > 0$ if the body
$K$ is a zonoid. A crucial ingredient in the proof of this result
is an embedding $\mathrm{K}_i$ of $\mathbf{Val}_i^+$ into the
space $C(\mathrm{Gr}_i)$ of continuous functions on the
Grassmannian $\mathrm{Gr}_i$ of $i$-dimensional subspaces of
$\mathbb{R}^n$ constructed by Klain \textbf{\cite{klain00}}:

Suppose that $\varphi \in \mathbf{Val}_i^+$, $1 \leq i \leq n -
1$. Then, by Lemma \ref{hadvol}, the restriction of $\varphi$ to
any subspace $E \in \mathrm{Gr}_i$ is proportional to the
$i$-dimensional volume $\mathrm{vol}_E$ on $E$, say
\[\varphi|_E = (\mathrm{K}_i\varphi)(E)\,\mathrm{vol}_E.  \]
The continuous function $\mathrm{K}_i\varphi : \mathrm{Gr}_i
\rightarrow \mathbb{R}$ defined in this way is called the {\it
Klain function} of $\varphi$. The induced map
\[\mathrm{K}_i: \mathbf{Val}_i^+ \rightarrow C(\mathrm{Gr}_i)  \]
is called the {\it Klain embedding}.

\begin{theorem} \label{klainthm} \emph{{\bf (Klain \cite{klain00})}}
The Klain embedding is injective.
\end{theorem}

Theorem \ref{klainthm} follows from a volume characterization of
Klain \textbf{\cite{klain95}}. Note, however, that the map
$\mathrm{K}_i$ is not onto; its image was described in terms of
the decomposition under the action of the group $\mathrm{SO}(n)$
by Alesker and Bernstein \textbf{\cite{AlBern}}.

The natural question how to reconstruct a valuation $\varphi \in
\mathbf{Val}_i^+$ given its Klain function $\mathrm{K}_i\varphi$,
was answered by Klain \textbf{\cite{klain00}} for centrally
symmetric convex sets. Since we need Klain's inversion formula
for zonoids only, we state just this special case. To this end,
we denote by $[u_1,\ldots,u_i]$ the $i$-dimensional volume of the
parallelotope spanned by $u_1, \ldots, u_i \in S^{n-1}$.

\begin{theorem} \label{invform} \emph{{\bf (Klain \cite{klain00})}}
Suppose that $\varphi \in \mathbf{Val}_i^+$ with $1 \leq i \leq n
- 1$. \linebreak If $Z_1\ldots,Z_i \in \mathcal{K}^n$ are zonoids
with generating measures $\mu_{Z_1}, \ldots, \mu_{Z_i}$, then
\[\varphi(Z_1,\ldots,Z_i)=\frac{1}{i!}\!\int_{S^{n-1}}\!\!\!\!\!\! \cdots\! \int_{S^{n-1}}
\!\!\! (\bar{\mathrm{K}}_i\varphi)(u_1,\ldots,u_i)
[u_1,\ldots,u_i] d\mu_{Z_1}(u_1)\cdots d\mu_{Z_i}(u_i),   \] where
\[(\bar{\mathrm{K}}_i\varphi)(u_1,\ldots,u_i) = \left \{ \begin{array}{ll}
(\mathrm{K}_i\varphi)(\mathrm{span}\{u_1, \ldots, u_i\}) &
\mbox{if } [u_1, \ldots, u_i] > 0, \\ 0 & \mbox{otherwise.}
\end{array} \right. \]
In particular, for any zonoid $Z \in \mathcal{K}^n$, we have
\[\varphi(Z)=\frac{1}{i!}\!\int_{S^{n-1}}\!\!\!\!\!\! \cdots\! \int_{S^{n-1}}
\!\!\!(\bar{\mathrm{K}}_i\varphi)(u_1,\ldots,u_i)\,
[u_1,\ldots,u_i]\,d\mu_Z(u_1)\cdots d\mu_Z(u_i).\]
\end{theorem}

\vspace{1cm}

\centerline{\large{\bf{ \setcounter{abschnitt}{4}
\arabic{abschnitt}. Proof of Theorem 2}}}

\reseteqn \alpheqn \setcounter{theorem}{0}

\vspace{0.6cm}

Before we can present the proof of Theorem 2, we need the
following auxiliary result.

\begin{lem} \label{equivconj} For $K \in \mathcal{K}^n$, the following statements are equivalent:
\begin{enumerate}
\item[(a)] For every non-negative $\varphi \in \mathbf{Val}$, its homogeneous components
$\varphi_i$ satisfy $\varphi_i(K) \geq 0$ for $0 \leq i \leq n$.
\item[(b)] For every $\Phi \in \mathbf{MVal}$, there exist $L_0, L_n \in \mathcal{K}^n$ (depending only on $\Phi$) and $\Phi_1(K), \ldots, \Phi_{n-1}(K) \in
\mathcal{K}^n$ such that (\ref{conjminkdecomp}) holds for every
$\lambda > 0$.
\end{enumerate}
\end{lem}

\noindent {\it Proof.} Let $K \in \mathcal{K}^n$ be fixed and
first assume that $\varphi_i(K) \geq 0$, $0 \leq i \leq n$, for
the homogeneous components $\varphi_i$ of any non-negative
$\varphi \in \mathbf{Val}$. Suppose that $\Phi \in
\mathbf{MVal}$. Then, by Lemma \ref{minkdecom}, for every $L \in
\mathcal{K}^n$, there are convex bodies $L_0, L_n \in
\mathcal{K}^n$ and continuous functions $g_i(L,\cdot)$ such that
\begin{equation} \label{equiv00}
h(\Phi(\lambda L),\cdot) = h(L_0,\cdot) +
\sum_{i=1}^{n-1}\lambda^i g_i(L,\cdot) + \lambda ^n
V(L)h(L_n,\cdot)
\end{equation}
for every $\lambda > 0$. In order to prove (b), it remains to show
that for each $i \in \{1, \ldots, n - 1\}$, the function
$g_i(K,\cdot)$ is the support function of a convex body
$\Phi_i(K)$. Since, by Lemma \ref{minkdecom}, the functions
$g_i(K,\cdot)$ are positively \linebreak homogeneous of degree
one, it suffices to prove that
\begin{equation} \label{equiv0}
g_i(K,x + y) \leq g_i(K,x) + g_i(K,y)
\end{equation}
for every $x,y \in \mathbb{R}^n$ and $i \in \{1, \ldots, n -
1\}$. To this end, fix $x,y \in \mathbb{R}^n$ and define $\psi \in
\mathbf{Val}$ by
\[\psi(L)=h(\Phi(L),x)+h(\Phi(L),y)-h(\Phi(L),x+y), \qquad L \in \mathcal{K}^n.  \]
Since support functions are sublinear, $\psi$ is non-negative.
Moreover, by (\ref{equiv00}), the homogeneous components
$\psi_i$, $1 \leq i \leq n - 1$, of $\psi$ are given by
\[\psi_i(L)=g_i(L,x)+g_i(L,y)-g_i(L,x+y).  \]
Since $\psi_i(K) \geq 0$ for $0 \leq i \leq n$, we obtain
(\ref{equiv0}). Thus, (a) implies (b).

\vspace{0.1cm}

Assume now that (b) holds. Suppose that $\varphi \in \mathbf{Val}$
is non-negative and let $\varphi_i$, $0 \leq i \leq n$ denote its
homogeneous components. Define a Minkowski valuation $\Phi \in
\mathbf{MVal}$ by
\[\Phi(L) = \varphi(L)B^n, \qquad L \in \mathcal{K}^n.  \]
Since $\varphi \geq 0$, the valuation $\Phi$ is well defined.
Using (\ref{homdecomp1}), it is easy to see that, on one hand,
\begin{equation} \label{equiv1}
h(\Phi(\lambda K),\cdot) = \varphi_0(K) + \lambda \varphi_1(K) +
\cdots + \lambda^n\varphi_n(K)
\end{equation}
for every $\lambda > 0$. On the other hand, it follows from (b)
that there exist $L_0, \Phi_1(K), \ldots, \Phi_{n-1}(K), L_n \in
\mathcal{K}^n$ such that
\begin{equation} \label{equiv2}
h(\Phi(\lambda K),\cdot) = h(L_0,\cdot) +
\sum_{i=1}^{n-1}\lambda^ih(\Phi_i(K),\cdot) +
\lambda^nV(K)h(L_n,\cdot)
\end{equation}
for every $\lambda > 0$. Comparing coefficients in (\ref{equiv1})
and (\ref{equiv2}) shows that \linebreak
$\varphi_0(K)=h(L_0,\cdot)$, $\varphi_n(K)=V(K)h(L_n,\cdot)$ and
$\varphi_i(K)=h(\Phi_i(K),\cdot)$ for $1 \leq i \leq n - 1$. This
is possible only if $\varphi_i(K) \geq 0$ for every $i \in \{0,
\ldots, n\}$. \hfill $\blacksquare$

\vspace{0.4cm}

Lemma \ref{equivconj} shows that Problem \ref{mainconj} is
equivalent to the question whether the homogeneous components of
any non-negative valuation in $\mathbf{Val}$ are also
non-negative. (A result of the latter type for monotone
valuations was \linebreak recently established by Bernig and Fu
\textbf{\cite{bernigfu10}}).

Using Theorem \ref{invform} and Lemma \ref{equivconj}, we now
establish an affirmative \linebreak answer to Problem
\ref{mainconj} for the class of zonoids:

\begin{theorem} \label{homdeczon} If $\Phi \in \mathbf{MVal}$,
then for every zonoid $Z \in \mathcal{K}^n$, there exist convex
bodies $L_0, \Phi_1(Z), \ldots, \Phi_{n-1}(Z), L_n \in
\mathcal{K}^n$ such that
\begin{equation*}
\Phi(\lambda Z) = L_0 + \lambda \Phi_1(Z) + \cdots +
\lambda^{n-1}\Phi_{n-1}(Z) + \lambda^nV(Z)L_n
\end{equation*}
for every $\lambda > 0$.

\end{theorem}

\noindent {\it Proof.} By Lemma \ref{equivconj}, it suffices to
show that $\varphi_i(Z) \geq 0$, $0 \leq i \leq n$, for the
homogeneous components $\varphi_i$ of any non-negative $\varphi
\in \mathbf{Val}$ and every zonoid $Z \in \mathcal{K}^n$. To this
end, first note that, by (\ref{homdecomp1}), for any $K \in
\mathcal{K}^n$,
\[0 \leq \varphi(\lambda K) = \varphi_0(K) + \lambda \varphi_1(K) + \cdots + \lambda^n\varphi_n(K)  \]
for every $\lambda > 0$. By letting $\lambda$ tend to zero, we
therefore, see that $\varphi_0$ is always non-negative for any
non-negative $\varphi \in \mathbf{Val}$. Similarly, dividing by
$\lambda^n$ and letting $\lambda$ tend to infinity, it follows
that $\varphi_n$ is always non-negative.

It remains to show that $\varphi_i(Z)\geq 0$, $1 \leq i \leq n -
1$, for any zonoid $Z \in \mathcal{K}^n$. In order to see this,
let $K \in \mathcal{K}^n$ be a centrally symmetric convex body
contained in an $i$-dimensional subspace $E$ with
$\mathrm{vol}_E(K) > 0$. By Lemma \ref{hadvol}, \linebreak we have
$\psi(K) = 0$ for any $\psi \in \mathbf{Val}_j$ with $j
> i$. Therefore, it follows that for any non-negative $\varphi \in \mathbf{Val}$,
\[0 \leq \varphi(\lambda K) = \varphi_0(K) + \lambda \varphi_1(K) + \cdots + \lambda^{i-1} \varphi_{i-1}(K) + \lambda^i \varphi_i(K)  \]
for every $\lambda > 0$. Again, dividing by $\lambda^i$ and
letting $\lambda$ tend to infinity, we see that $\varphi_i(K) \geq
0$. Let $\varphi_i^{\pm}$ denote the even and odd parts of
$\varphi_i$, respectively. Since $K$ is centrally symmetric, we
conclude $\varphi_i^-(K)=0$ and
\[0 \leq \varphi_i(K) = \varphi_i^+(K)= (\mathrm{K}_i\varphi_i^+)(E)\,\mathrm{vol}_E(K).  \]
Since the subspace $E$ was arbitrary, we see that
$\mathrm{K}_i\varphi_i^+ \geq 0$. Consequently, by Theorem
\ref{invform}, $\varphi_i^+(Z) \geq 0$ for any zonoid $Z \in
\mathcal{K}^n$. Moreover, since zonoids are centrally symmetric,
we have $\varphi_i^-(Z)=0$, and thus
$\varphi_i(Z)=\varphi_i^+(Z)\geq 0$. \phantom{aaa} \hfill
$\blacksquare$

\vspace{0.3cm}

Theorem 2 is now a simple consequence of Theorem \ref{homdeczon}.
It is the special case $Z = B^n$ of the following

\begin{koro} \label{minksteinzon} Suppose that $\Phi \in \mathbf{MVal}$ and let $K \in \mathcal{K}^n$.
Then for every zonoid $Z \in \mathcal{K}^n$ there exist (unique)
$\Phi_Z^{(j)} \in \mathbf{MVal}$ such that
\begin{equation} \label{minkzon}
 \Phi(K + rZ) = \sum_{j=0}^n r^{n-j} \Phi_Z^{(j)}(K)
\end{equation}
for every $r > 0$.
\end{koro}

\noindent {\it Proof.} Let $K \in \mathcal{K}^n$ be fixed and
define $\Psi^K: \mathcal{K}^n \rightarrow \mathcal{K}^n$ by
\[\Psi^K(L)=\Phi(K + L), \qquad L \in \mathcal{K}^n.  \]
It is easy to see that, in fact, $\Psi^K \in \mathbf{MVal}$. Thus,
by Theorem \ref{homdeczon}, for every zonoid $Z$, there exist
$\Psi_0^K(Z), \ldots, \Psi_n^K(Z) \in \mathcal{K}^n$ such that
\[\Phi(K + rZ)=\Psi^K(rZ) =  \Psi_0^K(Z)+r\Psi_1^K(Z)+\cdots + r^{n-1}\Psi_{n-1}^K(Z) + r^n\Psi_n^K(Z) \]
for every $r > 0$. Define $\Phi_Z^{(j)}: \mathcal{K}^n \rightarrow
\mathcal{K}^n$ by
\[\Phi_Z^{(j)}(L)=\Psi_{n-j}^L(Z), \qquad L \in \mathcal{K}^n.  \]
Clearly, the maps $\Phi_Z^{(j)}$ satisfy (\ref{minkzon}).
Moreover, from an application of the Steiner formula
(\ref{mcmullenstein}) to the valuation
$\varphi(K)=h(\Phi(K),\cdot)$ (with values in the vector space
$C(S^{n-1})$) and the uniqueness of the derived valuations
$\varphi^{(j)}$, it follows that $\Phi_Z^{(j)} \in \mathbf{MVal}$.
\hfill $\blacksquare$

\vspace{0.4cm}

We end this section with a further generalization of Theorem
\ref{homdeczon}.

\begin{theorem} \label{minkpoly} Suppose that $\Phi \in \mathbf{MVal}$ and let $Z_1, \ldots, Z_m \in \mathcal{K}^n$ be
zonoids. Then
\[\Phi(\lambda_1Z_1 + \cdots + \lambda_mZ_m), \qquad \lambda_1, \ldots, \lambda_m \geq 0,\]
can be expressed as a polynomial in $\lambda_1, \ldots, \lambda_m$
of total degree at most $n$ whose coefficients are convex bodies.

\end{theorem}

\noindent {\it Proof.} Using arguments as in the proof of Lemma
\ref{equivconj}, we see that it is enough to prove that
$\varphi_i(Z_{j_1},\ldots,Z_{j_i}) \geq 0$, $1 \leq j_1, \ldots,
j_i \leq m$, holds for the mixed valuations derived from any
non-negative valuation $\varphi_i \in \mathbf{Val}_i$ with $1
\leq i \leq n - 1$. To this end, let $\varphi_i^{\pm}$ denote the
even and odd parts of $\varphi_i$, respectively. In the proof of
Theorem \ref{homdeczon}, we have seen that
$\mathrm{K}_i\varphi_i^+ \geq 0$. Consequently, by Theorem
\ref{invform}, $\varphi_i^+(Z_{j_1},\ldots,Z_{j_i})\geq 0$.
Moreover, since zonoids are centrally symmetric, we have
$\varphi_i^-(Z_{j_1},\ldots,Z_{j_i})=0$. Thus, we conclude
$\varphi_i(Z_{j_1},\ldots,Z_{j_i})=\varphi_i^+(Z_{j_1},\ldots,Z_{j_i})\geq
0. $ \hfill $\blacksquare$

\vspace{1cm}

\centerline{\large{\bf{ \setcounter{abschnitt}{5}
\arabic{abschnitt}. Background material for the proof of Theorem 3
}}}

\reseteqn \alpheqn \setcounter{theorem}{0}

\vspace{0.6cm}

For quick reference, we state in the following the geometric
inequalities (for which we refer the reader to the book by
Schneider \textbf{\cite{schneider93}}) and other ingredients
needed in the proof of Theorem 3.

For $K, L \in \mathcal{K}^n$ and $0 \leq i \leq n$, we write
$V(K[i],L[n-i])$ for the mixed volume with $i$ copies of $K$ and
$n - i$ copies of $L$. For $K, K_1, \ldots, K_i \in \mathcal{K}^n$
\linebreak and $\mathbf{C}=(K_1,\ldots,K_i)$, we denote by
$V_i(K,\mathbf{C})$ the mixed volume
$V(K,\ldots,K,K_1,\ldots,K_i)$ where $K$ appears $n - i$ times.
For $0 \leq i \leq n - 1$, we use $W_i(K,L)$ to denote the mixed
volume $V(K[n-i-1],B^n[i],L)$. The mixed volume $W_i(K,K)$ will be
written as $W_i(K)$ and is called the \emph{$i$th \linebreak
quermassintegral} of $K$. The \emph{$i$th intrinsic volume}
$V_i(K)$ of $K$ is defined by
\begin{equation} \label{viwi}
\kappa_{n-i}V_i(K)=\binom{n}{i} W_{n-i}(K),
\end{equation}
where $\kappa_m$ is the $m$-dimensional volume of the Euclidean
unit ball in $\mathbb{R}^m$. A special case of (\ref{mixed}) is
the {\it classical Steiner formula} for the volume of the
parallel set of $K$ at distance $r > 0$,
\[V(K + rB^n)=\sum \limits_{i=0}^n r^i{n \choose i}W_i(K) = \sum
\limits_{i=0}^n r^{n-i}\kappa_{n-i}V_i(K).\]

Associated with a convex body $K \in \mathcal{K}^n$ is a family
of Borel measures $S_i(K,\cdot)$, $0 \leq i \leq n - 1$, on
$S^{n-1}$, called the {\it area measures of order} $i$ of $K$,
such that for each $L \in \mathcal{K}^n$,
\begin{equation} \label{defsi}
W_{n-1-i}(K,L)=\frac{1}{n}\int_{S^{n-1}} h(L,u)\,dS_i(K,u).
\end{equation}
In fact, the measure $S_i(K,\cdot)$ is uniquely determined by the
property that (\ref{defsi}) holds for all $L \in \mathcal{K}^n$.
The existence of a polynomial expansion of the translation
invariant valuation $W_{n-1-i}(\cdot,L)$, thus carries over to
the surface area measures. In particular, for $r > 0$, we have the
Steiner type formula
\begin{equation} \label{steinsi}
S_j(K + rB^n,\cdot)=\sum_{i=0}^j r^{j-i} {j \choose i}
S_{i}(K,\cdot).
\end{equation}

Let $\mathcal{K}^n_{\mathrm{o}}$ denote the set of convex bodies
in $\mathbb{R}^n$ with non-empty interior. \linebreak One of the
fundamental inequalities for mixed volumes is the general
Minkowski inequality: If $K, L \in \mathcal{K}^n_{\mathrm{o}}$
and $0 \leq i \leq n - 2$, then
\begin{equation} \label{genmink}
W_i(K,L)^{n-i} \geq W_i(K)^{n-i-1}W_i(L),
\end{equation}
with equality if and only if $K$ and $L$ are homothetic.

A consequence of the Minkowski inequality (\ref{genmink}) is the
Brunn--Minkowski inequality for quermassintegrals: If $K, L \in
\mathcal{K}^n_{\mathrm{o}}$ and $0 \leq i \leq n - 2$, then
\begin{equation} \label{quermassbm}
W_i(K+L)^{1/(n-i)} \geq W_i(K)^{1/(n-i)}+W_i(L)^{1/(n-i)},
\end{equation}
with equality if and only if $K$ and $L$ are homothetic.

A further generalization of inequality (\ref{quermassbm}) (where
the equality conditions are not yet known) is the following: If
$0 \leq i \leq n-2$, $K, L, K_1, \ldots, K_i \in \mathcal{K}^n$
and $\mathbf{C}=(K_1,...,K_i)$, then
\begin{equation} \label{mostgenbm}
V_i(K+L,\mathbf{C})^{1/(n-i)} \geq
V_i(K,\mathbf{C})^{1/(n-i)}+V_i(L,\mathbf{C})^{1/(n-i)}.
\end{equation}

The {\it Steiner point} $s(K)$ of $K \in \mathcal{K}^n$ is the
point in $K$ defined by
\[s(K)=n\int_{S^{n-1}}h(K,u)u\,du,  \]
where the integration is with respect to the rotation invariant
probability measure on $S^{n-1}$. It is not hard to show that $s$
is continuous, rigid motion equivariant and Minkowski additive.

\begin{theorem} \label{steinerpoint}
{\rm (see e.g., \textbf{\cite[\textnormal{p.\
307}]{schneider93}})} The Steiner point map $s: \mathcal K^n
\rightarrow \mathbb R^n$ is the unique vector valued rigid motion
equivariant and continuous valuation.
\end{theorem}

A convex body $K$ is said to be of class $C^2_+$ if the boundary
of $K$ is a $C^2$ submanifold of $\mathbb{R}^n$ with everywhere
positive Gaussian curvature. The following property of bodies of
class $C^2_+$ will be useful (see \textbf{\cite[\textnormal{p.\
150}]{schneider93}}):

\begin{lem} \label{c2plus} If $K \in \mathcal{K}^n$ is a convex body of class
$C_+^2$, then there exist a convex body $L \in \mathcal{K}^n$ and
a real number $r > 0$ such that $K = L  + rB^n$.
\end{lem}

In the remaining part of this section we recall the notion of
smooth \linebreak valuations as well as the definition of an
important derivation operator on the space $\mathbf{Val}$ needed
in the next section. In order to do this, we first endow the space
$\mathbf{Val}$ with the norm
\[\|\varphi \| = \sup \{|\varphi(K)|: K \subseteq B^n \}, \qquad \varphi \in \mathbf{Val}.  \]
It is well known that $\mathbf{Val}$ thus becomes a Banach space.
The group $\mathrm{GL}(n)$ acts on $\mathbf{Val}$ continuously by
\[(A \varphi)(K)=\varphi(A^{-1}K), \qquad A \in \mathrm{GL}(n),\, \varphi \in \mathbf{Val}. \]
Note that the subspaces $\mathbf{Val}_i^{\pm} \subseteq
\mathbf{Val}$ are invariant under this $\mathrm{GL}(n)$ action.
Actually they are also irreducible. This important fact was
established by Alesker \textbf{\cite{Alesker01}}; it directly
implies a conjecture by McMullen that the linear combinations of
mixed volumes are dense in $\mathbf{Val}$. A different dense subset of
$\mathbf{Val}$ can be defined as follows:

\vspace{0.3cm}

\noindent {\bf Definition} \emph{A valuation $\varphi \in \mathbf{Val}$ is called smooth if
the map $\mathrm{GL}(n) \rightarrow \mathbf{Val}$ defined by $A
\mapsto A\varphi$ is infinitely differentiable.}

\vspace{0.3cm}

We use $\mathbf{Val}^{\infty}$ to denote the space of continuous
translation invariant and smooth valuations. For the subspaces of
homogeneous valuations of given parity in $\mathbf{Val}^{\infty}$
we write $\mathbf{Val}_i^{\pm,\infty}$. It is well known (cf.\
\textbf{\cite[\textnormal{p.\ 32}]{wallach1}}) that
$\mathbf{Val}_i^{\pm,\infty}$ is a dense $\mathrm{GL}(n)$
invariant subspace of $\mathbf{Val}_i^{\pm}$. Moreover, from
Corollary \ref{mcmullen} one deduces that the space
$\mathbf{Val}^{\infty}$ admits a direct
 sum decomposition into its subspaces
of homogeneous valuations of given parity.

The Steiner formula (\ref{mcmullenstein}) gives rise to an
important derivation operator $\Lambda: \mathbf{Val} \rightarrow
\mathbf{Val}$ defined by
\[(\Lambda \varphi)(K) = \left . \frac{d}{dt} \right |_{t = 0} \varphi(K + t B^n).    \]

Note that $\Lambda$ commutes with the action of the orthogonal
group $\mathrm{O}(n)$. \linebreak Moreover, if $\varphi \in
\mathbf{Val}_i$ then $\Lambda \varphi \in \mathbf{Val}_{i-1}$.

The significance of the linear operator $\Lambda$ can be seen from
the following Hard Lefschetz type theorem established for even
valuations by Alesker \textbf{\cite{Alesker03}} and for general
valuations by Bernig and Br\"ocker \textbf{\cite{Bernig07b}}:

\begin{theorem} \label{hardlef} If $2i \geq n$, then $\Lambda^{2i-n}:
\mathbf{Val}^{\infty}_{i} \rightarrow
\mathbf{Val}^{\infty}_{n-i}$ is an isomorphism. \linebreak In
particular, $\Lambda: \mathbf{Val}^{\infty}_{i} \rightarrow
\mathbf{Val}^{\infty}_{i-1}$ is injective if $2i-1 \geq n$ and
surjective if $2i - 1 \leq n$.
\end{theorem}

Theorem \ref{hardlef} suggests that it is (sometimes) possible to
deduce results on valuations of degree $i$ from result on
valuations of a degree $j > i$. We will exploit this idea in the
proof of Theorem \ref{main2}. To this end, we use
Theorem~\ref{main1} \linebreak to define the derivation operator
$\Lambda$ for translation invariant continuous Minkowski
valuations:

\begin{koro} \label{derminkval}  Suppose that $\Phi \in \mathbf{MVal}$. Then there exists a $\Lambda
\Phi \in \mathbf{MVal}$ such that for every $K \in \mathcal{K}^n$
and $u \in S^{n-1}$,
\[h((\Lambda \Phi)(K),u) = \left . \frac{d}{dt} \right |_{t=0} h(\Phi(K+tB^n),u).   \]
Moreover, if $\Phi \in \mathbf{MVal}_i$ with $1 \leq i \leq n$,
then $\Lambda \Phi \in \mathbf{MVal}_{i-1}$.
\end{koro}

We remark that a Hard Lefschetz type theorem similar to Theorem
\ref{hardlef} does {\it not} hold for the sets $\mathbf{MVal}_i$
in general. This follows from results of Kiderlen
\textbf{\cite{kiderlen05}} and the second-named author
\textbf{\cite{Schu06a}} on translation invariant and
$\mathrm{SO}(n)$ equivariant Minkowski valuations of degree $1$
and $n - 1$, respectively.

\vspace{1cm}

\centerline{\large{\bf{ \setcounter{abschnitt}{6}
\arabic{abschnitt}. The symmetry of bivaluations}}}

\reseteqn \alpheqn \setcounter{theorem}{0}

\vspace{0.6cm}

We recall here the notion of bivaluations and, in particular, a recent result on the
symmetry of rigid motion invariant homogeneous bivaluations. As a reference
for the material in this section we refer to the recent article
\textbf{\cite{ABS2011}}.

\vspace{0.4cm}

\noindent {\bf Definition} A map $\phi: \mathcal{K}^n \times
\mathcal{K}^n \rightarrow \mathbb{R}$ is called a
\emph{bivaluation} if $\phi$ is a valuation in each argument. A
bivaluation $\phi$ is called {\it translation biinvariant} if
$\phi$ is invariant under independent translations of its
arguments and $\phi$ is said to be \emph{$\mathrm{O}(n)$
invariant} if $\phi(\vartheta K,\vartheta L) = \varphi(K,L)$ for
all $K, L \in \mathcal{K}^n$ and $\vartheta \in \mathrm{O}(n)$.
We say $\phi$ has {\it bidegree} $(i,j)$ if $\phi(\alpha K,\beta
L)=\alpha^i \beta^j\phi(K,L)$ for all $K,L \in \mathcal{K}^n$ and
$\alpha, \beta > 0$.

\vspace{0.4cm}

Important examples of invariant bivaluations can be constructed
using mixed volumes and Minkowski valuations:

\pagebreak

\noindent {\bf Examples:}

\begin{enumerate}
\item[(a)] For $0 \leq i \leq n$, the bivaluation $\phi: \mathcal{K}^n
\times \mathcal{K}^n \rightarrow \mathbb{R}$, defined by
\[\phi(K,L)=V(K[i],L[n-i]), \qquad K, L \in \mathcal{K}^n,  \]
is translation biinvariant and $\mathrm{O}(n)$ invariant and has
bidegree $(i,n-i)$.

\item[(b)] Suppose that $0 \leq i, j \leq n$ and let $\Phi_j \in \mathbf{MVal}_j$ be $\mathrm{O}(n)$ {\it
equivariant}, i.e.\ $\Phi_j(\vartheta K)= \vartheta \Phi_j(K)$ for
all $K \in \mathcal{K}^n$ and $\vartheta \in \mathrm{O}(n)$. Then
the bivaluation $\psi: \mathcal{K}^n \times \mathcal{K}^n
\rightarrow \mathbb{R}$, defined by
\[\psi(K,L)=W_{n-i-1}(K,\Phi_j(L)), \qquad K, L \in \mathcal{K}^n,   \]
is translation biinvariant and $\mathrm{O}(n)$ invariant and has
bidegree $(i,j)$.
\end{enumerate}

First classification results for invariant bivaluations were
obtained only recently by Ludwig \textbf{\cite{Ludwig10a}}.
Systematic investigations of continuous translation biinvariant
bivaluations were started in \textbf{\cite{ABS2011}}. In order to
describe some of the results obtained there, let $\mathbf{BVal}$
denote the vector space of all continuous \linebreak translation
biinvariant (real valued) bivaluations. We write
$\mathbf{BVal}_{i,j}$ for its subspace of all bivaluations of
bidegree $(i,j)$ and we use $\mathbf{BVal}^{\mathrm{O}(n)}$ and
$\mathbf{BVal}_{i,j}^{\mathrm{O}(n)}$ to denote the respective
subspaces of $\mathrm{O}(n)$ invariant bivaluations.

\vspace{0.15cm}

From McMullen's decomposition of the space $\mathbf{Val}$ stated
in Corollary \ref{mcmullen}, one immediately deduces a
corresponding result for the space $\mathbf{BVal}$:

\begin{koro} \label{mcmullenbval}
\[\mathbf{BVal}=\bigoplus_{i,j=0}^n \mathbf{BVal}_{i,j}.\]
\end{koro}

Minkowski valuations arise naturally from data about projections
or sections of convex bodies. For example, the projection body
operator \linebreak $\Pi \in \mathbf{MVal}_{n-1}$ is defined as
follows: The {\it projection body} $\Pi(K)$ of $K$ is the convex
body defined by $h(\Pi(K),u)=\mathrm{vol}_{n-1}(K|u^{\bot})$, $u
\in S^{n-1}$, where $K|u^{\bot}$ denotes the projection of $K$
onto the hyperplane orthogonal to $u$. Introduced already by
Minkowski, projection bodies have become an important tool in
several areas over the last decades (see, e.g.\
\textbf{\cite{habschu09, lutwak85, lutwak93, LYZ2000a,
schneider93, schnweil}}; for their special role in the theory of
valuations see \textbf{\cite{ludwig02, Ludwig10a, SchuWan11}}).

An extremely useful and well known symmetry property of the
projection body operator is the following: If $K, L \in
\mathcal{K}^n$, then
\begin{equation} \label{pisym}
V(\Pi(K),L,\ldots,L) = V(\Pi(L),K,\ldots,K).
\end{equation}
Variants of this identity and its generalizations have been used
extensively for establishing geometric inequalities related to
convex and star body valued valuations (see e.g.\
\textbf{\cite{gardner2ed, haberl08, habschu09, Ludwig10,
lutwak85, lutwak86, lutwak88, lutwak90, lutwak93, LYZ2000a,
LYZ200b, LZ1997,Schu06b, Schu06a, Schu09}}).

In \textbf{\cite{ABS2011}} the following generalization of the
symmetry property (\ref{pisym}) was established:

\vspace{0.1cm}

\begin{theorem} \label{thm_symmetry} If $\phi \in \mathbf{BVal}_{i,i}^{\mathrm{O}(n)}$, $0
\leq i \leq n$, then
\begin{equation*}
 \phi(K,L)=\phi(L,K)
\end{equation*}
for every $K, L \in \mathcal{K}^n$.
\end{theorem}

For $m = 1,2$ let the partial derivation operators $\Lambda_m:
\mathbf{BVal} \rightarrow \mathbf{BVal}$ be defined by (cf. Corollary \ref{derminkval})
\begin{equation} \label{deflam1}
(\Lambda_1 \phi)(K,L) = \left . \frac{d}{dt} \right |_{t = 0}
\phi(K + t B^n,L).
\end{equation}
and
\begin{equation} \label{deflam2}
(\Lambda_2 \phi)(K,L) = \left . \frac{d}{dt} \right |_{t = 0}
\phi(K,L + t B^n).
\end{equation}

\vspace{0.1cm}

\noindent Clearly, if $\phi \in \mathbf{BVal}_{i,j}$, then
$\Lambda_1 \phi \in \mathbf{BVal}_{i-1,j}$ and $\Lambda_2 \phi \in
\mathbf{BVal}_{i,j-1}$.

Also define an operator $\mathrm{T}: \mathbf{BVal} \rightarrow
\mathbf{BVal} $ by
\[(\mathrm{T}\phi)(K,L)=\phi(L,K).  \]
By Theorem \ref{thm_symmetry}, the restriction of $\mathrm{T}$ to
$\mathbf{BVal}_{i,i}^{\mathrm{O}(n)}$ acts as the identity. Thus,
we obtain the following immediate consequence of Theorem
\ref{thm_symmetry}:

\begin{koro} \label{comdiag} Suppose that $0 \leq j \leq n$ and $0\leq i \leq j$. Then the following diagram is commutative:
\[
\begin{xy}
 \xymatrix{
  \mathbf{BVal}_{j,j}^{\mathrm{O}(n)} \,\, \ar[dd]_{\textnormal{\normalsize $\Lambda_2^{j-i}$}} \ar[rr]^{\textnormal{\normalsize $\mathrm{T} = \mathrm{Id}$}}    &     & \,\,  \mathbf{BVal}_{j,j}^{\mathrm{O}(n)} \ar[dd]^{\textnormal{\normalsize $\Lambda_1^{j-i}$}}  \\
       &     &     \\
      \mathbf{BVal}_{j,i}^{\mathrm{O}(n)} \,\,   \ar[rr]^{\textnormal{\normalsize $\mathrm{T}$}} &     & \,\, \mathbf{BVal}_{i,j}^{\mathrm{O}(n)}
  }
\end{xy}
\]
\end{koro}

\noindent {\it Proof.} Suppose that $\phi \in
\mathbf{BVal}_{j,j}^{\mathrm{O}(n)}$, $0 \leq j \leq n$ and let
$K, L \in \mathcal{K}^n$. Then, by Theorem \ref{thm_symmetry}, we
have
\[\phi(L,K+tB^n) = \phi(K + tB^n,L)  \]
for every $t > 0$. Consequently, by definitions (\ref{deflam1})
and (\ref{deflam2}), we obtain
\[(\mathrm{T}\Lambda_2^{j-i}\phi)(K,L)=(\Lambda_1^{j-i} \phi)(K,L).   \]

\vspace{-0.7cm}

\hfill $\blacksquare$

\vspace{0.4cm}

In the proof of Theorem \ref{main2} we will repeatedly make
critical use of the following consequence of Corollary
\ref{comdiag}:

\begin{koro} \label{durch} Let $\Phi_j \in \mathbf{MVal}_j$, $2 \leq j \leq n - 1$, be $\mathrm{SO}(n)$ equivariant. If $1 \leq i \leq j + 1$, then
\[W_{n-i}(K,\Phi_j(L))=\frac{(i-1)!}{j!}W_{n-1-j}(L,(\Lambda^{j+1-i}\Phi_j)(K))  \]
for every $K, L \in \mathcal{K}^n$.
\end{koro}

\noindent {\it Proof.} We first note that any $\mathrm{SO}(n)$
equivariant $\Phi \in \mathbf{MVal}$ is also $\mathrm{O}(n)$
equivariant (see \textbf{\cite[\textnormal{Lemma
7.1}]{ABS2011}}). Define $\phi \in
\mathbf{BVal}_{j,j}^{\mathrm{O}(n)}$ by
\[\phi(K,L) = W_{n-1-j}(K,\Phi_j(L)), \qquad K, L \in \mathcal{K}^n.  \]
From (\ref{defsi}) and (\ref{steinsi}), it follows that
\[W_{n-i}(K,\Phi_j(L))=\frac{(i-1)!}{j!}(\Lambda_1^{j+1-i}\phi)(K,L).  \]
Thus, applications of Corollary \ref{comdiag} and Corollary
\ref{derminkval} complete the proof. \hfill $\blacksquare$

\vspace{1cm}

\centerline{\large{\bf{ \setcounter{abschnitt}{7}
\arabic{abschnitt}. Brunn--Minkowski type inequalities}}}

\reseteqn \alpheqn \setcounter{theorem}{0}

\vspace{0.6cm}

In this final section we present the proof of Theorem
\ref{main2}. Special cases of this result for $j = n - 1$ were
obtained in \textbf{\cite{Schu06b}} and for $j \in \{2, \ldots, n
- 1\}$ and $i = j + 1$ in \textbf{\cite{ABS2011}} and
\textbf{\cite{Schu09}}. The equality conditions for bodies of
class $C_+^2$ are new for $j \leq n - 2$.

\begin{theorem} \label{satzgenbmbm} Let $\Phi_j \in \mathbf{MVal}_j$, $2 \leq j \leq n - 1$, be $\mathrm{SO}(n)$
equivariant and non-trivial. If $K, L \in \mathcal{K}^n$ and $1
\leq i \leq j + 1$, then
\begin{equation} \label{desinequ}
W_{n-i}(\Phi_j(K+L))^{1/ij} \geq
W_{n-i}(\Phi_j(K))^{1/ij}+W_{n-i}(\Phi_j(L))^{1/ij}.
\end{equation}
If $K$ and $L$ are of class $C^2_+$, then equality holds if and
only if $K$ and $L$ are homothetic.
\end{theorem}
{\it Proof.} By Corollary  \ref{durch}, we have for $Q \in
\mathcal{K}^n$,
\begin{equation} \label{bm1}
W_{n-i}(Q,\Phi_j(K+L))=\frac{(i-1)!}{j!}W_{n-1-j}(K+L,(\Lambda^{j+1-i}\Phi_j)(Q)).
\end{equation}
From an application of inequality (\ref{mostgenbm}), we obtain
\begin{equation} \label{bm2}
\begin{array}{c}
W_{n-1-j}(K+L,(\Lambda^{j+1-i}\Phi_j)(Q))^{1/j} \qquad \qquad \qquad \qquad \mbox{\phantom{aaaaaaaaaaaaa}} \\
\phantom{aa} \geq
W_{n-1-j}(K,(\Lambda^{j+1-i}\Phi_j)(Q))^{1/j}+W_{n-1-j}(L,(\Lambda^{j+1-i}\Phi_j)(Q))^{1/j}.
\end{array}
\end{equation}
A combination of (\ref{bm1}) and (\ref{bm2}) and another
application of Corollary \ref{durch}, therefore show that
\begin{eqnarray} \label{cons1}
W_{n-i}(Q,\Phi_j(K+L))^{1/j} \geq  W_{n-i}(Q,\Phi_j(K))^{1/j}+
W_{n-i}(Q,\Phi_j(L))^{1/j}.
\end{eqnarray}
It follows from Minkowski's inequality (\ref{genmink}), that
\begin{equation} \label{mink1}
W_{n-i}(Q,\Phi_j(K))^i\geq W_{n-i}(Q)^{i-1}W_{n-i}(\Phi_j(K)),
\end{equation}
and, similarly,
\begin{equation} \label{mink2}
W_{n-i}(Q,\Phi_j(L))^i\geq W_{n-i}(Q)^{i-1}W_{n-i}(\Phi_j(L)).
\end{equation}
Plugging (\ref{mink1}) and (\ref{mink2}) into (\ref{cons1}), and
putting $Q=\Phi_j(K+L)$, now yields the desired inequality
(\ref{desinequ}).

In order to derive the equality conditions for convex bodies of
class $C^2_+$, we first show that such bodies are mapped by
$\Phi_j$ to bodies with non-empty interior. For the following
argument, the authors are obliged to T.~Wannerer. Let $Q \in
\mathcal{K}^n$ be of class $C^2_+$. By Lemma \ref{c2plus}, there
exist a convex body $M \in \mathcal{K}^n$ and a number $r > 0$
such that $Q = M + rB^n$. Using that $\Phi_j$ has degree $j$, we
thus obtain from Theorem \ref{main1} that
\[\Phi_j(Q)= \Phi_j(M + rB^n)= r^j \Phi_j^{(0)}(M) + \cdots + \Phi_j^{(j)}(M),  \]
where $\Phi_j^{(m)} \in \mathbf{MVal}_m$, $1 \leq m \leq j$, and
$\Phi_j^{(0)}(M)=\Phi_j(B^n)$. Since $\Phi_j$ is $\mathrm{SO}(n)$
equivariant, we must have $\Phi_j(B^n)=tB^n$ for some $t \geq 0$.
Since $\Phi_j$ is non-trivial, it follows from an application of
Corollary \ref{durch} (with $i = 1$ and $K = B^n$) that in fact
$t > 0$. Consequently, $\Phi_j(Q) \in \mathcal{K}_{\mathrm{o}}^n$.

Now assume that $K, L \in \mathcal{K}^n$ are of class $C^2_+$ and
that equality holds in (\ref{desinequ}). Let $s$ be the Steiner point map.
Since $Q \mapsto s(\Phi_j(Q)) + s(Q)$ is a continuous and rigid motion equivariant valuation,
Theorem \ref{steinerpoint} implies that $s(\Phi_j(Q)) = 0$ for every $Q \in \mathcal{K}^n$.
Thus, we deduce from the equality conditions of (\ref{mink1}) and
(\ref{mink2}), that there exist $\lambda_1, \lambda_2 > 0$ such
that
\begin{equation} \label{hom2}
\Phi_j(K)=\lambda_1\Phi_j(K+L) \qquad \mbox{and} \qquad
\Phi_jL=\lambda_2 \Phi_j(K+L)
\end{equation}
and
\[\lambda_1^{1/j}+\lambda_2^{1/j} = 1.   \]
Moreover, (\ref{hom2}) and another application of Corollary
\ref{durch} (with $i = 1$ and $K = B^n$) imply that
\[W_{n-j}(K)=\lambda_1W_{n-j}(K+L) \qquad \mbox{and} \qquad  W_{n-j}(L)=\lambda_2W_{n-j}(K+L).  \]
Hence, we have
\[W_{n-j}(K+L)^{1/j}=W_{n-j}(K)^{1/j}+W_{n-j}(L)^{1/j},  \]
which, by (\ref{quermassbm}), implies that $K$ and $L$ are
homothetic. \hfill $\blacksquare$

\vspace{0.3cm}

We remark that our proof shows that the smoothness assumption in
the equality conditions can be omitted for bodies with non-empty
interior in case $\Phi_j$ maps those bodies again to convex
bodies with non-empty interior. This is always the case if $j = n
- 1$ (cf.\ \textbf{\cite{Schu06b}}), but is not known in general.

We also note that since translation invariant continuous Minkowski
\linebreak valuations of degree one are linear with respect to
Minkowski addition (see e.g.\ \textbf{\cite{hadwiger51}}),
inequality (\ref{desinequ}) also holds in the case $j=1$ (this
follows from (\ref{quermassbm})). The equality conditions,
however, are different in this case.

\vspace{0.5cm}

\noindent {{\bf Acknowledgments} The work of the authors was
supported by the \linebreak Austrian  Science Fund (FWF), within
the project ``Minkowski valuations and geometric inequalities",
Project number: P\,22388-N13.

\begin{small}

Vienna University of Technology \par Institute of Discrete
Mathematics and Geometry \par Wiedner Hauptstra\ss e 8--10/1046
\par A--1040 Vienna, Austria


\end{small}

\end{document}